# Paradox regained: Life beyond Gödel's shadow

*Bhupinder Singh Anand*

*Gödel's explicit thesis was that his undecidable formula* **GUS** *is a well-formed, well-defined proposition in any formalisation of Intuitive Arithmetic* **IA** *in which the Axioms and Rules of Inference are recursively definable. His implicit thesis was that* **GUS** *is not formally inconsistent in any such system.*

*I consider a constructive, and intuitionistically unobjectionable, formalisation* **PP** *of* **IA** *in which the Axioms and Rules of Inference are recursively definable. I argue that, although Gödel's proposition* **GUS** *is a well-formed formula in* **PP**, *it is an ill-defined proposition that formally reflects the "Liar" paradox in* **PP**.

*I argue that the introduction of "formal truth" and "**PP**-provability" values to selected propositions of* **IA** *through* **PP** *leads to the collapse of the Gödelian argument advanced by J.R.Lucas, Roger Penrose and others.*

## CONTENTS

### Author's Preface

### 1. Introduction

## 6. Questioning Gödel's conclusions

  6.1. *GDL*(*x*) is well-defined in **PA**, but not necessarily in any interpretation **M**
  6.2. **GUS** is well-formed and well-defined, but Gödel's reasoning is non-constructive
  6.3. *Intuitive Arithmetic* **IA** is not a model of **PA**
  6.4. *Generalisation* as a *Rule of Inference* is intuitionistically questionable
  6.5. A non-constructive interpretation of *Generalisation*

## 7. Conclusions

## Author's preface

*Intuitive truth - a paradigm shift*

The alternative consequences of Godel's reasoning developed in this paper (*and in its roots*) involve a significant paradigm shift in the philosophical basis of our perception of the nature of *Intuitive Knowledge*, and more particularly of the concept of factual, or intuitive truth.

By *Intuitive Knowledge* I refer loosely to that body of pro-active knowledge that stems directly from our conscious states, in contrast to our reactive *Instinctive Knowledge*, which stems from, and lies within, our sub-conscious and unconscious states.

This preface is intended to highlight the wider significance of an issue that may otherwise lie in obscurity due to the specialised nature of the subject. My thesis is that the influence, on our current modes of thought, of the interpretations, and conclusions, drawn from Gödel's original paper (*Gödel 1931*) may have a wider, multi-disciplinary, element that is not obvious from an appreciation of its purely logical and mathematical import.

*Implicit influence of Gödel's Platonism*

I argue that the roots of this influence may be traced to implicitly Platonistic elements that underlie classical first order predicate calculus **PA**, which is based essentially on the formal system defined by Gödel in his (*1931*) paper. Loosely speaking Gödel, who was an explicitly strong Platonist, assumed the existence of a world of ideals that could objectively be referred to for arbitrating which of our assumptions or premises, when formally expressed in a rigorously constructed scientific language, were "intuitively true", and which were "intuitively false", under interpretation.

Now one may, of course, argue reasonably - as Gödel does - in Platonistic terms and define "intuitive truths" as characteristics of "relationships" that are assumed to "exist" in some "absolute" sense (*that is, even in the absence of any "perceiver"*) between the "objects" of an external "ontology" (*both of which are also taken to "exist" in some "absolute" sense*).

I argue, however, for an alternative view of "relationships" as belonging to "perceptions" that we consciously "construct" and selectively "assign" to abstract "objects" (*that themselves are conceptual "constructs"*) of an abstract "ontology" (*that is similarly a conceptual "construct"*).

In other words, I argue that each "perception" can reasonably be assumed to be an abstract "construct" based on a unique, one-of-a-kind, never-to-be-repeated consciousness of an "experience". "Intuitive truth" is then, essentially, a "constructed", space-time-localised, "factual truth". It corresponds to a "subjectively" constructed characteristic of the "expression" of individually constructed "perceptions". Loosely speaking, it corresponds to a characteristic of the way we construct an "expression" for that which we select as common to a series of "subjective perceptions", rather than to a characteristic that we "discover" of an "objectively" observed "something".

My brief against Platonism, when rooted in scientific ways of thought, is that it could permit us to "logically" validate our subjective "intuitive" perceptions as being reflective of some "absolute Truth" that "must be" of universal "significance" in a Utopian, Platonistic world.

*Gödel's Theorems and "formal (logical) truth"*

The view of Gödel's Theorems that I attempt to present is that they are essentially concerned with the effectiveness of our ability to "communicate" the abstractions that we intellectually "construct" on the basis of our "individualistic" sensory perceptions. They thus have to do with the efficiency and effectiveness of our language of communication, and involve the concept of "logical truth", or "formal truth" (*the two terms are used synonymously in this paper*).

I argue that, in formal languages, a selected set of "axiomatic truths" is expressed as a set of *Axioms* (*or Axiom schemas*). The selection criteria is that the *Axioms* are readily accepted by any "perceiver" as faithfully reflecting some "significant" "factual truths" pertaining to the expression of abstract "constructed" elements of the "constructed" ontology under consideration as "perceived" and "conceived" also by the "perceiver".

For the most basic, and intuitive, of our scientific languages, namely *Number Theory* or our *Intuitive Arithmetic* of the natural numbers, I take the commonly accepted selection of such "axiomatic truths" as the set of *Peano's Postulates*, first expressed in semi-axiomatic format by Dedekind (*1901*).

I consider that the challenge in any such theory is then to find a suitable set of *Rules of Inference* by which we can assign unique "formal (*logical*) truth" values to as many well-formed propositions of the theory as possible that are not *Axioms*.

I argue that the concept of "formal (*logical*) truths" is thus merely the result of the application of a set of *Rules of Inference* for assigning such "formal (*logical*) truth" values to various "logical" permutations and combinations of "axiomatic truths" as (*finite* *and infinite*) compound assertions (*which ideally should not introduce any new "axiomatic truths" that are not already implicit within the Axioms*).

I argue that the *Rules of Inference* of any specific language should include not only the familiar "*Rules of Inference*" that we define for the language, but also the logical *Axioms* that, strictly speaking, should be considered equally as being part of the *Rules of Inference* of the language.

*Significance*

I argue that he significance of "formal (*logical*) truths" lies in our experience that the "factual (*intuitive*) truths" of our "perceptions" can generally be corresponded in a communicable language with a high degree of correlation to the "formal (*logical*) truths" of the language. The large body of our "factual (*intuitive*) truths" may thus also be viewed as "intuitively constructed" permutations and combinations of a smaller, ideally finite, set of "factual (*intuitive*) truths".

This appears to suggest that the "significance" to us (*and possibly to any intelligence whose evolution is based on communication*) of any set of "factual (*intuitive*) truths" may be proportional to the body of "formal (*logical*) truths" that can be inferred by various "logical" combinations and permutations of the "axiomatic truths" that can be distilled from the particular set of "factual (*intuitive*) truths".

I argue that the significance of Gödel's Theorems lies in the fact that they are derived in a system of *Axioms* where the *Rules of Inference* lead to a particularly rich body of expressions that can be assigned "formal (*logical*) truth" values under various interpretations of the symbols of the theory.

*Non-verifiable assertions implicit in Gödel's formal system*

I argue, however, that a major feature of Gödel's formal system is that the chosen *Rules of Inference* are non-constructive, in the sense that they also assign, implicitly and sweepingly, non-verifiable "formal (*logical*) truth" values under interpretation in some models to various expressions. Thus the system, in a sense, postulates "formally (*logically*) true" expressions under some interpretation that cannot be correlated, even in principle, to any "factual (*intuitive*) truths" of a human "perception".

In (*Anand 2001*), I highlight this non-intuitive aspect of Gödel's choice of *Rules of Inference*, and suggest more intuitive, constructive, *Rules of Inference* that would make our formal systems of *Arithmetic* more faithfully representative of our *Intuitive Arithmetic*.

One of the consequences of using constructive *Rules of Inference* is that the reasoning by which Godel concludes that the formal "Gödel" proposition, which translates loosely as "The 'Gödel' proposition is not 'provable'", is intuitively true collapses. In this paper, I define a formalisation of *Intuitive Arithmetic* in which Gödel's argument leads to a contradiction. In effect, we are faced with the expected reflection of the "Liar" sentence, namely "The 'Liar' sentence is a lie", in constructive formal theories.

However, I argue that this is not necessarily the drawback that it appears to be at first sight. In fact such contradictions both encourage and discipline us in our use of rich and creative languages. They force us to focus on devising rules by which we can recognise well-formed formulas of the theory that need to be treated as ill-defined propositions (*in the sense that there is no way a unique "formal (logical) truth" value can be assigned to the well-formed formula by the Rules of Inference of the theory*). I thus offer the thesis that "consistency" may be treated in any theory as the specification for determining which well-formed formulas qualify as well-defined propositions.

*An illusion*

I argue that, in a sense, what Gödel's reasoning actually accomplishes is best viewed as an "illusion", where a "constructive" argument successfully masks a "non-constructive" premise.

The significance of this is that if an "illusion-prone" logical system forms the bedrock of all scientific thought and discourse, then attendant conflicts and controversies must inevitably follow "illusory", if not outright "hallucinatory", "perceptions" rooted in non-constructive thought.

To the extent that the foundations of our scientific thought can visibly be seen to be increasingly influencing, and impacting on, the way we "perceive", "conceive", "express", "communicate" and "act" in every sphere of human activity, any logic that provides a theoretical legitimacy to "illusory perceptions" as having a possible validity in a non-verifiable universe that "exists" in some "absolute" sense can be held qualitatively accountable in some measure for reinforcing the kind of "perception" that leads to unresolvable confrontations arising out of "extra-logical" consequences of non-constructive theories.

(*The development of systems of Artificial Intelligence, which so far have not shown any marked susceptibility to Gödelian "illusions", may, however, force consideration of formal systems that are more self-disciplined in their faithfulness to our intuitively constructive aspects that can be verified and duplicated mechanically.*)

*To sum up*

In this paper I aim to highlight the implicit assumptions that, I argue, underlie the formal structure of our current formal system of logic.

*Truth*: I argue, albeit implicitly, against the concept of "intuitive truth" as something a priori that is to be discovered, and which can be captured only partially through formal structures. In other words, I argue against the concept that "provable" formulas in any formal system can be corresponded 1-1 to a proper subset of the "intuitively true" propositions under any "interpretation" of the formal system.

*Formal truth*: My thesis is that "formal (*logical*) truth" is a characteristic that we choose to attribute, as something of commonly communicable, intellectual, "significance", to our assertions through meta-rules of assignment that we call our logical *Rules of Inference* (*by definition, these lie outside the language*).

*Significance*: I hold we choose these rules in order to intellectually communicate the essence of "factual (*intuitive*) truths" that we also choose to intuitively construct as "significant", individualistic, abstract concepts.

*Factual truth*: I hold that "factual truths", or "intuitive truths", belong to our consciously constructed intuitive "perceptions", just as "formal (*logical*) truth" belongs to our consciously constructed intellectual systems (*languages*) of communication.

*Consciousness*: I hold that "factual (*intuitive*) truth" does not "exist" outside "consciousness". Thus there is no "factual (*intuitive*) truth" in the absence of conscious intuitive "perception", just as there is no "formal (*logical*) truth" in the absence of conscious intellectual "communication".

*Axiomatic truth*: I hold that "axiomatic truths" are "factual (*intuitive*) truths" that we choose to express within a language of communication in order to create a common link between the body of our intuitively constructed "factual (*intuitive*) truths", and the body of our intellectually constructed "formal (*logical*) truths".

*Rules of Inference*: I hold that our *Rules of Inference* and *Logical Axioms* are what we select as the means of intellectually constructing the body of our "formal (*logical*) truths", just as our *Intuition* is what we rely on for intuitively constructing the body of our "factual (*intuitive*) truths".

*Inductive truth*: I further hold that "inductive (*extra-logical*) truth" is, loosely speaking, an individualistic, non-experiential, non-verifiable, extrapolation of "formal (*logical*) truth" - built around both qualified, and individually preferred, generalisations of "formal (*logical*) truth".

*Interpretation*: I argue against the thesis that a formal system can be meaningfully viewed in isolation, separately from an intended interpretation. I argue that, given a set of "axiomatic truths", it is our *Rules of Inference* that determine "formal (*logical*) truth" values in any interpretation.

*Model*: Thus I argue (*Anand 2001, Section_4*) that an interpretation can be a model of two formal systems that have identical primitive symbols, identical rules of construction for well-formed formulas, and identical sets of "axiomatic truths", but distinctly different, contradictory, "logical" consequences.

*Communication*: My thesis is that it is through our system of formal languages that we, on the one hand, constructively assign "formal (*logical*) truths" (*and provability*) and, on the other, non-constructively assign selected characteristics of non-verifiable "inductive (*extra-logical*) truths" (*through qualified generalisation*), to elements of the language in an effort to correspond these to specific "factual (*intuitive*) truths" of our intuitively constructed individual perceptions whose "significance" we desire to communicate.

*Uniqueness*: I suggest that there can be no commonly verifiable non-standard model. Each such must differ in some essentials perceived uniquely only by each perceiver who is, in effect, both the architect and builder of the model.

*Bhupinder Singh Anand*                                        *Friday, 15 February 2002*

# Paradox regained: Life beyond Gödel's shadow

## 1. Introduction

### 1.1. An intellectual challenge for the 21st century: Gödel in a wider perspective

One of the intellectual challenges inherited from the last century is that of interpreting Gödel's assertions of formal undecidability within a wider, intuitively acceptable, perspective. Given his formal premises, Gödel's formal reasoning and its formal consequences (*Gödel 1931*) are logically irrefutable. However, there is no similarly unequivocal interpretation of his reasoning, or of his conclusions.

I argue elsewhere (*Anand 2001*), in greater detail, the thesis that differing perceptions of the interpretation of Gödel's reasoning and conclusions arise because the system of *standard Peano's Arithmetic*, essentially formulated by Gödel (*1931, pp.9*), is only one - and perhaps not the most representative - of several significantly differing systems that can be defined to formalise our system of *Intuitive Arithmetic* **IA** of the natural "counting" numbers (*which we take to be the Arithmetic based on an intuitive interpretation of Dedekind's formulation of Peano's Postulates*[1]).

(*Another thesis would be that Gödel's formal system itself is open to significantly differing interpretations.*)

---

[1] These can be formulated (*Mendelson 1964*) as follows:

(**P1**)   0 is a natural number
(**P2**)   If $x$ is a natural number, there is another natural number, denoted by $x'$ (called the successor of $x$).
(**P3**)   0 is not equal to $x'$ for any natural number $x$.
(**P4**)   If $x' = y'$, then $x = y$.
(**P5**)   If **Q** is a property which may or may not hold of natural numbers, and if (I) 0 has the property **Q**, and (II) whenever a natural number $x$ has the property **Q**, then $x'$ has the property **Q**, then all natural numbers have the property **Q** (*Principle of Induction*).

In this paper I outline in more general terms a constructive, and intuitionistically unobjectionable, formalisation **PP** of **IA** in which the *Axioms* and *Rules of Inference* are recursively definable. I argue that, although Gödel's undecidable sentence **GUS** is a well-formed formula in **PP**, it is an ill-defined formal proposition[2] that formally reflects the "Liar" paradox in **PP**.

**1.2. Gödel's Undecidability Theorem and its consequences**

Usually referred to loosely as "Gödel's Theorems", these assertions are the outcome of Gödel's attempt to construct a formal language for faithfully expressing our *Intuitive Arithmetic* **IA** of the natural numbers.

The ideal source for assessing Gödel's intent and the significance of the wealth of concepts introduced by him in the course of his attempt is still, albeit arguably, his seminal paper "On formally undecidable propositions of Principia Mathematica and related systems I", accessible on the web in an English translation by B. Meltzer (*Gödel 1931*).

However, with the advantage of hindsight, my thesis in this paper is that the essence of Gödel's reasoning, its formal consequences, and its wider implications, are better viewed today in a somewhat different context.

**1.3. The "Liar" paradox**

Around 2000 years back, Greek philosophers discussed the paradox of ambiguous self-reference within ordinary language. This arises if we postulate a 'Liar' expression, which reads as "The 'Liar' proposition is a lie", as a valid proposition within the language in which the expression is constructed. Clearly, the "Liar" proposition is true if and only if it is a lie!

---

[2] We define a "proposition", or a "sentence", as a well-formed formula of a formal system that contains no free variables.

Now the grammar of an ordinary language governs the construction of words and propositions of the language for use as a means of general unrestricted communication. It obviously does not contain the specification necessary to prohibit the creation by definition of such vague, or ill-defined, self-referential expressions as that of the "Liar". Though these may have the formal form of propositions, they appear devoid of the "communicative content", or "meaning", that we intuitively expect propositions to have.

We conclude that the well-defined "Liar" expression cannot be a well-defined "Liar" proposition in a "consistent" language. So we may either treat the language as inherently "inconsistent", or conclude that the language is deficient in its ability to identify well-defined expressions that are ill-defined, or "meaningless", propositions.

Considered primarily a linguistic anomaly, the paradox does, however, focus attention on the issue:

> When may we treat a well-formed expression as a well-defined proposition, without fear of contradiction?

## 1.4. The Russell paradox

The "Liar" paradox gained in significance when, in 1901, Bertrand Russell discovered a similarly paradoxical expression within Set Theory. Loosely speaking, he defined a set, which we may name "Russell", as the set whose elements are all, and only, those sets of the Theory that do not belong to themselves. Clearly, if it can be expressed within the Theory, then the set "Russell" belongs to itself if and only if it does not belong to itself!

Again, the *Axioms* and logical *Rules of Inference* of the Theory - which governed the construction of sets for use in a more restricted and precise language of mathematical communication - did not contain the specification necessary to satisfactorily prohibit the creation by definition of vague, ill-defined, or inconsistent entities such as the "Russell" set that could be expressed within the language of the Theory, but apparently harboured concealed concepts involving unruly infinite sets.

Here too, we conclude from the contradiction that introduction of the well-defined "Russell" expression cannot yield a well-defined "Russell" set in a consistent Set Theory. So, again, we may either treat any Theory admitting such sets as inherently inconsistent, or conclude that such a Theory is deficient in its ability to identify well-defined expressions that yield ill-defined "sets".

Viewed as a reflection on the soundness of the foundations of mathematics, Russell's paradox focuses on the issue:

> When may we introduce entities through definition into our mathematical languages, without fear of contradiction?

**1.5. Creation through definition**

We note that definitions are essentially arbitrary assignments of convenient names within a language, or theory, for expressing complex reasoning in a compact and easier-to-grasp form.

The paradoxes indicate that unrestricted definitions, particularly those that involve self-reference, may admit well-formed expressions within the language, or theory, which lead to an inconsistency. Clearly such expressions must be considered as ill-defined elements within the language, or theory. It follows that a language, or theory, that admits such well-formed expressions as representing well-defined elements cannot uncritically be assumed consistent.

## 2. Formal systems

### 2.1. *Intuitive Arithmetic* **IA** and the formal system **PP**

A natural question arises whether the simplest, and most intuitive, of our mathematical theories also admit such expressions. We consider, for instance, the *Intuitive Arithmetic* **IA** of the natural "counting" numbers. The following *Axioms* - based on what are termed as *Peano's Postulates*, and expressed formally in a sub-language **PP** of **IA** as defined below - are commonly taken to be "formally true" representations of some of the most "intuitively true" assertions of **IA**[3]:

(*By "PP", we shall henceforth mean "the sub-language **PP** of **IA**".*)

(***PP1***)   $(Ax)|=_{PP} \sim(0 = (x+1))$

   (*Intuitive interpretation: Adding 1 to any natural number expressible in **PP** never yields 0.*)

---

[3] However, we note that, in the usual formalisations of Peano Arithmetics, such as **PP**, the axioms do not completely reflect Dedekind's intentions as formulated in his Peano Postulates. Thus there are no axioms corresponding to the postulates (*P1*) and (*P2*) detailed in footnote (1). These essentially assert that the domain of the system, over which the variables of the system range, necessarily consists of only the natural numbers. A corresponding axiom in a formalisation would assert that "$(Ax)('x$ is a numeral')". In **PP**, such an assertion could be expressed as "$(Ax)((x = 0) \mathbf{v} (E!y)(E!z)((x = (y+1)) \& (z = (x+1))))$".

It can reasonably be argued that it is the omission of this restricting axiom that leads to the questionable - and controversial - admission of non-constructive elements, and of non-standard interpretations, into what are intended essentially to be constructive and intuitionistically unobjectionable formalisations of our *Intuitive Arithmetic* of the natural "counting" numbers, as expressed by Dedekind in his Peano Postulates.

(**PP2**)  $(Ax)(Ay)\models_{PP} \sim(x=y) \Rightarrow \sim((x+1) = (y+1))$

   (*Intuitive interpretation: Adding 1 to each of two different natural numbers expressible in **PP** yields two different natural numbers that are expressible in **PP**.*)

(**PP3**)  $(Ax)\models_{PP} (x+0) = x$

   (*Intuitive interpretation: Adding 0 to a natural number that is expressible in **PP** yields the same natural number.*)

(**PP4**)  $(Ax)(Ay)\models_{PP} (x+(y+1)) = ((x+y)+1)$

   (*Intuitive interpretation: Adding 1 to a natural number that is expressible in **PP**, and then adding another natural number that is expressible in **PP** to the sum, yields the same natural number that is expressible in **PP** as is obtained by adding the two natural number that are expressible in **PP** first, and then adding 1 to their sum; in other words addition is an "associative" operation over the natural numbers expressible in **PP**.*)

(**PP5**)  $(Ax)\models_{PP} (x*0) = 0$

   (*Intuitive interpretation: Multiplying a natural number expressible in **PP** by 0 yields 0.*)

(**PP6**)  $(Ax)(Ay)\models_{PP} (x*(y+1)) = ((x*y)+x)$

   (*Intuitive interpretation: Adding 1 to a natural number that is expressible in **PP**, and then multiplying the sum by a second natural number that is expressible in **PP**, yields the same natural number that is expressible in **PP** as is obtained by multiplying the two natural number that are expressible in*

> *PP first, and then adding the second natural number that is expressible in PP to their product; in other words multiplication is a "distributive" operation over "addition" for the natural numbers expressible in PP.)*

## 2.2. The alphabet S of the sub-language PP of IA

The characteristic feature of these selected "assertions", which are taken to express the most intuitive arithmetical truths of **IA** symbolically, is that they are expressed using only a limited alphabet **S**, containing the arithmetical symbols "+" (*addition*) and "*" (*multiplication*) only, apart from the meta-logical and logical symbols "|=$_{PP}$" (*It is "formally true" in **PP** that*), "~"(*not*), "=>" (*implies*), "=" (*equals*), "&" (*and*), "**v**" (*or*), "(A$x$)" (*for all x*), "(E$x$)" (*there exists x*), "$x, y, ...$" (*variables*), "a, b, c, ..." (*constants*), "0" (*zero*), "1" (*one*), and the two parentheses "(", ")".

## 2.3. The domains of IA and PP

The (*ontological*) domain over which the finite set of variables "$x, y, ...$" range is taken, in **IA**, to be the intuitively non-terminating series "0, 1, 2, 3, ..." of natural "counting" numbers and, in **PP**, the non-terminating formal series expressing the natural numbers in **PP** as the numerals "0, 0+1, (0+1)+1, ((0+1)+1)+1, ...".

## 2.4. Quantification and "formal truth" in PP : nailing an ambiguity

Thus each of (*PP1*) to (*PP6*) symbolically assigns a "formal (*logical*) truth" value in **PP** to a countable infinity of arithmetical "propositions" of **IA**. In other words, assuming the expression, including the meta-logical symbol "|=$_{PP}$", is a well-formed "propositional formula" by the rules for the formation of well-formed "propositional formulas" of **PP**, the well-formed formula to the right of the meta-logical symbol "|=$_{PP}$" is asserted as a "formally true" proposition in **PP** for any given set of finite values of the variables contained in it.

Clearly, a crucial characteristic of the language **PP**, as defined above, is that all formulas containing quantifiers necessarily contain the "formally true" meta-logical symbol "$\models_{PP}$" immediately after the quantifier. This reflects the intention that, in **PP**, the quantifiers "(A$x$)" and "(E$x$)" should indeed "quantify" the "formal truth" of the expression "*F*($x$)" over a specific domain in order to yield a well-defined expression that has the symbolic form of a "proposition". Thus all formal propositions of **PP** that contain quantifiers are meta-assertions under interpretation.

I argue further that "(A$x$)$\models_{PP}$ *F*($x$)" is intended to formally represent the intuitive meta-assertion that "*F*($x$) is 'formally true' in **PP** for all values of $x$ that are expressible in **PP**". I therefore adopt the expanded form to highlight the significance of the intended meta-assertion (*or negation*) of the compound "formal truth", involved in "(A$x$)$\models_{PP}$ *F*($x$)", when we intend to adjoin "(A$x$)" to "*F*($x$)" in **PP**.

## 2.5. *Rules of inference* in PP

The two *Rules of Inference* in **PP**, for deriving propositions in **PP** that are "logical" consequences of the above *Postulates*, are *Modus Ponens* and *Induction*.

(***PPR1***)  *Modus Ponens*: From "(A$x$)$\models_{PP}$ *F*($x$)" and "(A$x$)$\models_{PP}$ (*F*($x$) => *G*($x$))" we may conclude "(A$x$)$\models_{PP}$ *G*($x$)".

(*Intuitive interpretation: If **F**($x$) is "formally true" for all natural numbers $x$ expressible in **PP**, and **G**($x$) is "formally true" whenever **F**($x$) is "formally true" for any natural number $x$ expressible in **PP**, then **G**($x$) is "formally true" for all natural numbers $x$ expressible in **PP**.*)

(***PPR2***)  *Induction*: From "$\models_{PP}$ *F*(0)" and "(A$x$)$\models_{PP}$ (*F*($x$) => *F*($x$+1))" we may conclude "(A$x$)$\models_{PP}$ *F*($x$)".

(*Intuitive interpretation: If $F(0)$ is "formally true", and $F(x+1)$ is "formally true" whenever $F(x)$ is "formally true" for any natural number $x$ expressible in **PP**, then $F(x)$ is "formally true" for all natural numbers $x$ expressible in **PP**.*)

## 3. Gödel's reasoning

### 3.1. Is there a "Liar" expression in IA?

Around 1930, Gödel considered whether the "Liar" expression, "The 'Liar' proposition is a lie", would be "reflected" formally in **PP** by asserting an arithmetical "Gödel" proposition that is "equivalent", loosely speaking, to the assertion "The 'Gödel' proposition is not 'provable'".

What Gödel had noticed, with remarkable insight, was that terms such as "well-formed formula"[4], "well-formed proposition" and "proof sequence" could be formally defined in the language of **PP**. He therefore argued that he could express meta-assertions about **PP** such as " '$F(x)$' is a well-formed formula in **PP**", " '$(Ax)|=_{PP} F(x)$' is a well-formed proposition in **PP**", " '$F(x)$' is a 'provable' formula in **PP**", " '$(Ax)|=_{PP} F(x)$' is a 'provable' proposition in **PP**" and "The 'Gödel' proposition is not 'provable' in **PP**", amongst others, as equivalent arithmetical propositions in **PP**.

### 3.2. Provability in PP

The essence of Gödel's argument, when applied to **PP**, is that we can define a proposition $P$ of **IA** as "**PP**-provable" if and only if there is a finite "proof sequence" consisting of propositions of **PP** each of which is either a "formally true" axiom (***PP1-6***), or a "formally true" immediate consequence of the preceding "formally true" propositions by the *Rules of Inference* (***PPR1-2***), and where $P$ is the final proposition in the sequence.

---

[4] We use the formal term "formula" as corresponding to the intuitive term "expression".

The "**PP**-provable" propositions are thus characterised by the fact that they are all expressed as "well-formed propositions", such as "$(Ax)|=_{PP} F(x)$", using only a small set of primitive, undefined symbols of the alphabet **S**. Clearly, the *Axioms* (***PP1-6***) are all "**PP**-provable" propositions.

### 3.3. Recursive functions in PP

The significance of the alphabet **S** selected for expressing "**PP**-provable" propositions is that many significant arithmetical functions of **IA** such as "*n*!" (*factorial*), "*m^n*" (*exponential*), "*m/n*" (*division*), "*n* is a prime number", amongst others, are defined recursively[5].

Thus we note that "*n*!", for instance, is defined by "0!=1" and "(*n*+1)!=(*n*+1)*(*n*!) for all natural numbers *n*". If we attempt to eliminate the symbol "!" on the right side in this definition, we soon discover that the function is not directly reducible to a form that is expressible in only the symbols of the alphabet **S**.

Now an obvious way to express "*n*!" as a function in **PP** is, of course, to introduce "!" as an additional primitive symbol into the alphabet **S**. However, since we can define infinitely many recursive functions that are similarly not expressible directly in **S**, every introduction of a new function by definition through a new "symbol" such as "!" would

---

[5] Loosely speaking, we assume that, in an Arithmetical system such as **IA**, a function or relation containing "free" variables is "recursive" if and only if, for any given set of values for the "free" variables in the definition of the function or relation, the arithmetical value of the function, or the "formal truth/falsity" of the relation, can be determined in a finite number of steps from the *Axioms* of the system using its *Rules of Inference* by some mechanical procedure.

Another assumed characteristic of recursive expressions in an Arithmetical system such as **IA** is that, for any given set of values for the "free" variables in the definition of the expression, it can be reduced, in a finite number of steps by some mechanical procedure, to an expression that consists of only a finite number of primitive symbols of a suitably constructed sub-language of **IA** such as **PP**, even though the definition and expression of the function or relation may involve an element of self-reference.

require a fresh determination as to whether the enlarged system " '**PP**'+'**!**' " is capable of yielding a logical paradox.

### 3.4. A Representation lemma of Hilbert and Bernays

This issue is addressed by a Representation lemma of Hilbert and Bernays that is based on defining a suitable Gödel "Beta-function" for every recursive function $f(x)$ definable in **IA**. By means of this stratagem, we can then establish that every recursive function $f(x)$ definable in **IA** is "equivalent" to some formal arithmetical relation $F(x, y)$ of **IA** that can be expressed in **PP**, in the sense that, for any natural numbers $k, m$:

(*i*)     If $f(k)=m$ is "formally true" in **IA**, then "$\models_{PP} F(\underline{k}, \underline{m})$"[6] is "**PP**-provable", whilst

(*ii*)     If $f(k)=m$ is "formally false" in **IA**, then "$\models_{PP} {\sim}F(\underline{k}, \underline{m})$" is "**PP**-provable".

### 3.5. Are recursive functions "consistent"?

The question of whether "new" recursive functions are "consistent"[7], in other words whether they can introduce into **IA** self-contradicting expressions similar to the "Liar" expression, is thus reduced, in a sense, to the question of whether the well-formed formulas of **IA** that are "**PP**-provable" are "consistent".

This was among the various specific issues that Gödel addressed, albeit in a slightly different language **PA** (*known as standard Peano's Arithmetic*), which was essentially intended to reflect the intuitively true propositions of **IA** most faithfully. However, before considering the details of Gödel's reasoning in standard **PA**, we take advantage of hindsight to consider his reasoning in the system **PP**.

---

[6] We denote by "$\underline{n}$" the formal "numeral" corresponding to the intuitive natural number "*n*".

[7] Loosely speaking, we assume that a system is "consistent" if there is no proposition "*F*" such that both "*F*" and "${\sim}F$" are "formally true" in the system.

### 3.6. Gödel-numbering

A key concept underlying Gödel's reasoning arises from the rather ordinary fact that every expression '*F*' constructed by concatenation from the primitive, undefined symbols of **S** can be assigned a unique natural number, which we term as the "Gödel" number of the expression '*F*'. Gödel's extraordinary achievement was to recognise that this fact could be used to reflect various meta-mathematical statements about **PP** as arithmetical propositions within **PP** (*Gödel 1931, pp.13*).

For instance, Gödel constructively established (*Gödel 1931, pp.17-22*) that we can define a recursive *prf*($x$, $y$) in **IA** (*Gödel's 'yBx'*), constructed out of 44 "simpler" recursive, such that, for any natural numbers $k$ and $m$, *prf*($k$, $m$) is "formally true" if and only if $k$ is the Gödel-number of a finite proof sequence $K$ in **PP** for some well-formed proposition $M$ in **PP** whose Gödel-number is $m$.

(*This is generally expressed as the assertion that the Axioms and Rules of Inference of **PP** are recursively definable or recursively enumerable.*)

### 3.7. Gödel's Self-reference Lemma

From the recursive *prf*($x$, $y$), he could then define another recursive *q*($x$, $y$) in **IA** which is "formally true" if and only if $x$ is the Gödel-number of a well-formed formula $H(z)$ of **PP** with a single free variable $z$, and $y$ is the Gödel-number of a proof of $H(x)$ in **PP**.

Hence the constructive self-reference, which lies at the core of Gödel's reasoning, is that, given any natural numbers $k$, $m$:

> *q*($h$. $j$) is "formally true" in **IA** <=> The formula $J$ in **PP**, whose Gödel number is $j$, is
> a proof sequence in **PP** of the well-formed proposition $H(\underline{h})$ of **PP**, where $h$ is the Gödel-number of the formula $H(z)$ in **PP**.

### 3.8. Hilbert & Bernays Representation Lemma

Now, by Hilbert and Bernays Representation lemma, $q(x, y)$ is constructively equivalent in **IA** to a unique well-formed formula $Q(x, y)$ of **PP**, expressible in **S**, such that for every pair of natural numbers $k, m$:

(*i*)   $q(k, m)$ is "formally true" in **IA** => "$|=_{PP} Q(k, m)$" is "**PP**-provable"

(*ii*)   $q(k, m)$ is "formally false" in **IA** => "$|=_{PP} {\sim}Q(k, m)$" is "**PP**-provable"

### 3.9. The undecidable "Gödel" proposition in PP

If $p$ is the Gödel-number of the well-formed formula "$(Ay)|=_{PP}({\sim}Q(x, y))$", we consider then the well-formed "Gödel" proposition **GUS** expressed by "$(Ay)|=_{PP}({\sim}Q(p, y))$".

(*a*) We assume firstly that $r$ is the Gödel-number of some proof-sequence **R** in **PP** for the proposition "$(Ay)|=_{PP}({\sim}Q(p, y))$". By Gödel's Self-reference lemma, $q(p, r)$ is "formally true" in **IA**. Also, by the Representation Lemma, this implies that "$|=_{PP} Q(p, r)$" is "**PP**-provable". However, assuming standard logical axioms for **PP**, from the "**PP**-provability" of "$(Ay)|=_{PP}({\sim}Q(p, y))$" we have that "$|=_{PP} {\sim}Q(p, r)$" is "**PP**-provable". It follows that there is no natural number $r$ that is the Gödel-number of a proof-sequence **R** in **PP** for the proposition "$(Ay)|=_{PP}({\sim}Q(p, y))$". Hence "$(Ay)|=_{PP}({\sim}Q(p, y))$" is not "**PP**-provable".

(*b*) We assume next that $r$ is the Gödel-number of some proof-sequence **R** in **PP** for the proposition "${\sim}(Ay)|=_{PP}({\sim}Q(p, y))$". It then follows that "$(Ey)|=_{PP} Q(p, y)$" is "**PP**-provable". Hence, assuming standard logical axioms for **PP**, "$|=_{PP} Q(p, r)$" is "**PP**-provable" for some natural number $r$. However, we have by (*a*) that $q(p, r)$ is "formally false" in **IA** for all *r*, and so "$\models_{PP} \sim Q(p, r)$" is "**PP**-provable" for all *r*. The contradiction establishes that "$\sim(Ay)\models_{PP}(\sim Q(p, y))$" is not "**PP**-provable".

(*c*) We conclude that the "**PP**-provability" of the proposition "$(Ay)\models_{PP}(\sim Q(p, y))$" is "undecidable", since neither "$(Ay)\models_{PP}(\sim Q(p, y))$" nor "$\sim(Ay)\models_{PP}(\sim Q(p, y))$" is "**PP**-provable".

## 4. Beyond Gödel

### 4.1. "Formal truth" and "L-provability" as formally assigned values

Now, as noted in §2.4, given a well-constructed language **L** that expresses relations between the various terms of an ontology, one way of viewing *Axioms* and *Rules of Inference* is as the means by which we formally assign properties such as "formal truth" and "**L**-provability" to the various formulas of the language that are defined as well-formed propositions.

As adopted in this paper, this view is, in essence, a non-Platonistic approach to the concept of "intuitive truth". It explicitly holds "formal truth", or "logical truth", to be an assignment of values to specified formulas of **IA**, where the assignment, if intuitive, is necessarily axiomatic and, if formal, follows from the *Axioms* by the *Rules of Inference* in a constructive and intuitionistically unobjectionable manner.

### 4.2. Effectiveness of a language

A natural question, then, is whether a given set of *Axioms* and *Rules of Inference* of a language **L** suffice to unequivocally assign a unique "formal truth" value to every formula that is a well-formed proposition, a concept that we may define as *semantic effectiveness*.

Another question would be whether such *Axioms* and *Rules of Inference* of the language further suffice to constructively determine whether every well-formed proposition that is "formally true" is "**L**-provable", which we may define as *syntactic effectiveness*.

The "Liar" sentence and the "Russell" set establish that intuitive *Axioms* and *Rules of Inference* of our ordinary languages of communication and of set theory - whether implicit or explicit - do not suffice to ensure "*semantic effectiveness*" for these languages.

The "Gödel" proposition **GUS** establishes that the above part of our *Intuitive Arithmetic* **IA** of the natural numbers, which is formalised by the *Axioms* (***PP1-6***) and the *Rules of Inference* (***PPR1-2***), is not "*syntactically effective*".

**4.3. Collapse of the "Gödelian" argument**

Now, by §3.9(*a*), we have that "$|=_{PP} \sim Q(p, r)$" is "**PP**-provable" for all *r*. We thus have that "$\sim Q(p, y)$" is "formally true" in **IA** for all *y*, since the domain of **IA** is expressible in **PP**, and every "**PP**-provable" proposition is, by definition, "formally true" in **IA**. It now follows by our intuitive definition of quantification in **IA**, that "$(Ay)(\sim Q(p, y))$" is assertable as a "formally true" proposition in **IA**. (*since the domain of IA is expressible in PP*), whence it follows that "$(Ay)|=_{PP}(\sim Q(p, y))$" is a "formally true" proposition in **IA**.

We thus have that, though "$(Ay)|=_{PP}(\sim Q(p, y))$" is not "**PP**-provable", we yet have "$(Ay)|=_{PP}(\sim Q(p, y))$" established as a "formally true" proposition in **IA** by the *Axioms* (***PP1-6***) and *Rules of Inference* (***PPR1-2***) of **PP**. However, the "formal truth" of the assertion in **IA** is clearly of a "definitional" nature, and a formal consequence of the *Axioms* and *Rules of Inference* of **PP**.

We note that this is a curious "formal" feature of "Gödelian" propositions in **PP**, which essentially demolishes the Gödelian argument (*see also Anand 2001, §1.13*). As offered by J.R.Lucas, Roger Penrose and others, this argument is the thesis that there is some non-mechanistic element - knowledge of which is Platonistically available to human intelligence but cannot be reflected in any machine intelligence - that is involved in establishing that a well-formed formula such as "$(Ay)|=_{PP}(\sim Q(p, y))$" is not "**PP**-provable", yet translates into a "formally true" proposition under interpretation in **IA**, which is defined as the "standard" model of **PP**.

### 4.4. Constructive "PP-provability"

Now it seems natural to consider whether we can augment either the *Axioms* (*PP1-6*) or the *Rules of Inference* (*PPR1-2*) so that we can assign a "**PP**-provability" value to "$(Ay)|=_{PP}(\sim Q(p, y))$", and thereby obtain a larger set of "**PP**-provable" propositions of **IA**.

In the light of the above reasoning, a naturally obvious way to achieve this would be to introduce as a *Rule of Inference*:

> (***PPR3***)   *Constructive **PP**-provability*: From "$(Ax)[|=_{PP}$ '$x$ is a numeral' $=> |=_{PP} F(x)]$" we may conclude "$(Ax)|=_{PP} F(x)$".

(*We note that (**PPR3**) is recursively definable in **PP** since "x is a numeral" is recursively definable in **PP** by "(x=0) v (E!y)(E!z)((x = (y+1)) & (z = (x+1)))", where "E!" denotes uniqueness of the existential assertion. We also note that, in a more formal exposition, we would distinguish between "**PP**" and "**PP**+**PPR3**".*)

### 4.5. Paradox regained : "*Constructive PP-provability*" and "Gödelian" propositions

We now have, by §3.9 (*a*), that $q(p, r)$ is "formally false" in **IA** for all $r$, and so "$|=_{PP} \sim Q(p, r)$" is "**PP**-provable" for all $r$. It then follows from (***PPR3***) that "$(Ay)|=_{PP}(\sim Q(p, y))$" is "**PP**-provable", and so no longer undecidable.

We thus have that:

(*i*)     "(A*y*)|=$_{PP}$(~***Q***(*p*, *y*))" is "**PP**-provable" => "(A*y*)|=$_{PP}$(~***Q***(*p*, *y*))" is not "**PP**-provable", and

(*ii*)     "(A*y*)|=$_{PP}$(~***Q***(*p*, *y*))" is not "**PP**-provable" => "(A*y*)|=$_{PP}$(~***Q***(*p*, *y*))" is "**PP**-provable".

So we have succeeded in establishing the "Gödelian" proposition "(A*y*)|=$_{PP}$(~***Q***(*p*, *y*))" as reflecting the "Liar" sentence in a formal sub-language of **IA** that has a *Constructive Provability Rule of Inference*. We conclude that the well-formed "Gödelian" proposition "(A*y*)|=$_{PP}$(~***Q***(*p*, *y*))" is essentially an ill-defined proposition in such a system.

We note that the above reasoning remains logically valid even if we eliminate "|=$_{PP}$" as a primitive symbol in **PP**. However, this would again obscure my thesis that "quantification" essentially assigns "formal truth" values to the set of propositions of **IA** that are determined by the *Axioms* (***PP1-6***) and *Rules of Inference* (***PPR1-2***) of **PP**.

**4.6. Beyond Gödel's shadow**

We are now faced with the question of whether to accept the *Constructive **PP**-provability Rule of Inference* as intuitively natural to our *Intuitive Arithmetic* **IA** or not. In the first case, we are faced with the dilemma of a well-formed, but ill-defined, proposition in **IA**. In the second, we are faced with the dilemma of accepting a system of *Intuitive Arithmetic* **IA**, but rejecting its inconvenient intuitive elements.

This point assumes significance when we see how Gödel addressed this issue in the sub-language **PA** of **IA**.

(*By "**PA**", we shall henceforth mean "the sub-language **PA** of **IA**".*)

If the above argument, establishing the existence of Gödel's undecidable proposition **GUS** as a well-formed, but ill-defined, proposition in constructive and intuitionistically unobjectionable formalisations of *Peano's Postulates* such as **PP**, is substantive, then the question arises:

> What feature of Gödel's formalisation of *Peano's Postulates* permits him to conclude the existence of **GUS** as both a well-formed formula *and* as a well-defined formal proposition in standard **PA**?

## 5. Basis of Gödel's conclusions

### 5.1. *Intuitive Arithmetic* IA and Peano's Arithmetic PA

We start by noting that, in standard (*first order*) *Peano's Arithmetic* **PA** (*Mendelson 1964, pp. 102*), which is essentially the system considered by Gödel (*1931*), the *Axioms* and *Rules of Inference* assign "**PA**-provability" values to selected well-formed formulas of **IA**. The standard *Axioms* of **PA** are:

(***GP1***) $\vdash_{PA} \sim(0 = (x+1))$

(***GP2***) $\vdash_{PA} \sim(x = y) \Rightarrow \sim((x+1) = (y+1))$

(***GP3***) $\vdash_{PA} (x+0) = x$

(***GP4***) $\vdash_{PA} (x+(y+1)) = ((x+y)+1)$

(***GP5***) $\vdash_{PA} (x*0) = 0$

(***GP6***) $\vdash_{PA} (x*(y+1)) = ((x*y)+x)$

The genesis of the *Axioms* (***GP1-6***) can clearly be seen to lie in (***PP1-6***); the formulas are again expressed using only a small set of undefined, primitive, symbols such as "+" (*addition*) and "*" (*multiplication*) only apart from the meta-logical and logical symbols "|-$_{PA}$" (*It is provable in* **PA** *that*), "~"(*not*), "=>" (*implies*), "=" (*equals*), "&" (*and*), "**v**" (*or*), "(A*x*)" (*for all x*), "(E*x*)" (*there exists x*), "*x, y, ...*" (*variables*), "a, b, c, ..." (*constants*), "0" (*zero*), "1" (*one*), and the two parentheses "(", ")".

The non-terminating series of numerals "0, 0+1, (0+1)+1, ((0+1)+1)+1, ..." is again taken as formally representing in **PA** the various, intuitively non-terminating, natural "counting" number series such as "0, 1, 10, 11, ..." (*in binary format*), or "0, 1, 2, 3, ..." (*in the more common decimal format*).

### 5.2. *Rules of inference* in **PA**

Gödel's *Rules of Inference* in **PA**, for deriving other "**PA**-provable" formulas from the *Axioms* (***GP1-6***), are *Modus Pon*ens, *Induction* and *Generalisation*.

Like the *Axioms* of **PA**, the first two clearly reflect their roots in (***PPR1-2***).

- (***GPR1***)   *Modus ponens*: From "|-$_{PA}$ *F*" and "|-$_{PA}$ (*F* => *G*)" we may conclude "|-$_{PA}$ *G*", where *F* and *G* are any well-defined formulas of **PA**.

- (***GPR2***)   *Induction*: From "|-$_{PA}$ *F*(0)" and "|-$_{PA}$ (A*x*)(*F*(*x*) => *F*(*x*+1))" we may conclude "|-$_{PA}$ (A*x*)*F*(*x*)".

### 5.3. Provability in PA

Now we note that, since the domain[8] of the variables is not specified in **PA**, each of (***GP1***) to (***GP6***) formally represents an "indeterminate" set of "**PA**-provable" arithmetical formulas of **IA**, where:

---

[8] Compare §2.3.

(*a*) the expression to the right of the meta-logical symbol "|-$_{PA}$" is a well-defined "formula" of **PA**, formed from the alphabet **S**, that does not contain the meta-logical symbol "|-$_{PA}$"; in other words it is a validly defined formula consisting only of a finite number of the primitive logical and arithmetical symbols of **S**, and constructed by some well-defined rules of construction of **PA**.

(*b*) a "formula" is defined "**PA**-provable" if and only if there is a finite sequence of formulas of **PA** each of which is either an *Axiom* of **PA** or an immediate consequence of the preceding formulas by the *Rules of Inference* of **PA**.

## 5.4. Significant features of PP and PA

We note that an underlying thesis of this paper is that **PP** and **PA** are both sub-languages of our *Intuitive Arithmetic* **IA** that attempt to formally express our "intuitively true" assertions of **IA** as "formally true", "**PP**-provable" and "**PA**-provable" propositions in "platform-independent"[9] languages of communication. As also noted earlier, the semantical and syntactical "*effectiveness*" of a sub-language may be taken as a measure of its expressive strength.

We now note some significant features that differentiate **PP** from **PA** in this respect:

(*a*) **PP** explicitly assigns both "formal truth" values and formal "**PP**-provability" values to formulas of **IA** that are built entirely from the finite set of undefined primitive meta-logical, logical and arithmetical symbols of **PP**, provided these are well-formed simple or "compound" propositions as defined by the rules for proposition formation in **PP**. We note that since, as determined in §4.5, there are ill-defined Gödelian sentences that are "**PP**-provable", but cannot be assigned any "formal truth" value, the "**PP**-provable" formulas of **IA** are not a proper subset of the "formally true" formulas of **IA**.

---

[9] We consider any assignment of "formal truth" and "provability" values as "platform-independent" if the determination of such values can be arrived at by some mechanical procedure that is constructive and intuitionistically unobjectionable.

(*b*) **PA**, as defined so far, assigns only "**PA**-provability" values to formulas built entirely from the finite set of undefined primitive symbols of **PA**, provided these are well-formed formulas as defined by the rules for formation of formulas of **PA**.

(*c*) Since "**PA**-provability" of a formula does not require that it have the form of a "proposition", there are "**PA**-provable" formulas that are not "propositions". It follows that there are "**PA**-provable" formulas that cannot be assigned any "formal truth" values in **IA** (*and, ipso facto, in PP*) under interpretation, and so the "**PA**-provable" formulas too are not a proper sub-set of the "formally true" formulas of **IA** (*and, ipso facto, of PP*).

(*d*) The issue then is to see:

   (*i*) whether every "**PA**-provable" formula that has the form of a "proposition" can be assigned a unique "formal truth" value in **IA** (*and, ipso facto, in PP*) under interpretation and, if so,

   (*ii*) whether we can provide an intuitionistically unobjectionable rule for inferring from a "**PA**-provable" formula that does not have the form of a "proposition" some "**PA**-provable" formula that does, so that the latter can be assigned a unique "formal truth" value in **IA** (*and, ipso facto, in PP*).

(*e*) The particular *Rule of Inference* selected by Gödel to achieve this last is termed as *Generalisation*. It is undisputedly accepted today as a critical element of the system of standard **PA**.

## 5.5. *Generalisation*: Gödel's rule of "Extrapolation"

In its usual form, *Generalisation* is expressed as:

(**GPR3**)   *Generalisation*: From "$\vdash_{PA} F(x)$" we may conclude "$\vdash_{PA} (Ax)F(x)$".

However, I argue that Gödel's *Generalisation Rule of Inference* is essentially a non-constructive rule of "extrapolation". It is a crucial introduction of a concept that does not draw upon our intuition for its legitimacy. From the verifiable "**PA**-provability" of formulas that are not propositions, it asserts by extrapolation the "formal (*logical*) truth" of propositions under some interpretation which may not be constructively verifiable. In other words, it extends the boundaries of what is intuitively familiar into the area of what appears familiar, but is in effect unverifiable, by a "logical" extrapolation.

I argue that the non-constructive and Platonistic nature of *Generalisation* becomes apparent if we recast the above as:

(**GPR4**)   *Generalisation\**: From "$\vdash_{PA} F(x)$" we may conclude "$\vdash_{PA} ((Ax)|=_M F(x))$".

Now what *Generalisation\** essentially asserts is that from the formal "**PA**-provability" of the formula "$F(x)$", we may conclude that "$F(x)$" is "formally (*logically*) true" for all values of $x$ in *any* "interpretation"[10] **M** of **PA**. *Generalisation\** is thus the means of assigning "formal (*logical*) truth" values to selected formulas of **PA** not only in **IA** (*and, ipso facto, in PP, the formal sub-language of IA that formalises the concept of "intuitive truth" in IA*), but in *any* interpretation **M** of **PA**.

---

[10] Loosely speaking, an "interpretation" **M** of a formal system such as **PA** is a 1-1 mapping of the undefined symbols of **PA** into the symbols of **M** under which each predicate symbol of **PA** corresponds to a "similar" relation of **M**, each function symbol of **PA** corresponds to a "similar" function in **M**, and each constant symbol of **PA** corresponds to some fixed element of **M** (*Mendelson 1964, pp. 49*).

Before we consider the wider implications of a *Rule of Inference* such as *Generalisation*, we briefly review Gödel's argument for establishing an undecidable proposition in **PA**.

**5.6. Gödel's undecidable proposition GUS**

As detailed in §3.9, if $p$ is the Gödel-number of the formula "(A$y$)(~$Q$($x$, $y$))", then **GUS** is the proposition "(A$y$)(~$Q$($p$, $y$))".

Clearly, **GUS** is a well-constructed arithmetical formula, which is uniquely identified by a well-constructed natural number $g$ (*its Gödel-number*).

Gödel's argument (*essentially along the lines outlined in §3.6 –3.9*) is then that if $r$ is the Gödel-number of some proof-sequence ***R*** in **PA** of the formula "(A$y$)(~$Q$($p$, $y$))", this can constructively be shown to imply the "formal truth", in **IA** (*and, ipso facto, in PP*), of another arithmetical formula "*q*($p$, $r$)" of **IA** which, in turn, can constructively be shown equivalent to the following assertion ***G*** that can also be expressed as a formula in **IA**:

(***G***)  "A proof for the sentence whose Gödel number is $p$ cannot be written out in a finite number of steps from the *Axioms* (***GP1-6***), using only the *Rules of Inference* (***GPR1-3***)".

Assuming **PA** is consistent, Gödel could then conclude:

(1)  **GUS** is "**PA**-provable" => "*q*($p$, $r$)" is "formally true" in **IA** => ***G*** is "formally true" in **IA** => **GUS** is not "**PA**-provable",

Hence **GUS** is not "**PA**-provable";

(2)  **~GUS** is "**PA**-provable" => **GUS** is "**PA**-provable"

Hence **~GUS** is not "**PA**-provable".

### 5.7. The roots of GUS

Now the involved (*but, contrary to appearances, not critical*) part of Gödel's reasoning is that "$q(p, r)$" is the instantiation of an intuitive primitive recursive arithmetic relation "$q(p, y)$" where $p$ is the Gödel-number of the formula "$(Ay)(\sim Q(x, y))$".

Also, **GUS** is the formula "$(Ay)(\sim Q(p, y))$" where we have, by the application of the Representation Lemma §3.8 to **PA**, that, for all natural numbers $r$:

(*i*)   "$q(p, r)$" is "formally true" in **IA** => "$Q(p, r)$" is "**PA**-provable", and

(*ii*)   "$q(p, r)$" is "formally false" in **IA** => "$\sim Q(p, r)$" is "**PA**-provable".

As in the "Liar" and "Russell" cases, both the formulas "$q(p, y)$" and "$Q(p, y)$" are well-formed formulas within **IA** and **PA** respectively. Hence the other formulas such as "$(Ay)(\sim Q(p, y))$", or **GUS**, and "$q(p, r)$" are also well-formed formulas.

### 5.8. Well-formed, well-defined formulas of Gödel's formal system of standard PA

Although this is not immediately obvious, Gödel's entire chain of reasoning, which establishes **GUS** as undecidable, critically rests on proving firstly that the well-formed formula "$(Ay)(\sim Q(p, y))$" is also a well-defined proposition in **PA**, and secondly that it translates as a well-defined proposition in **IA**, which is taken to be the standard interpretation of standard **PA**.

Classically, the well-definedness of "$(Ay)(\sim Q(p, y))$" is established by showing that "$(E!w)Q(x, w)$" is a "**PA**-provable" formula (*Mendelson 1964, pp. 134*).

(*Here "E!" denotes uniqueness of the existential assertion. The above argument is developed in detail in Anand 2001 (§2.8)*)

## 6. Questioning Gödel's conclusions

### 6.1. GUS is well-defined in PA, but not necessarily in any interpretation M

The questionable aspect of Gödel's reasoning and conclusions is that, since "$Q(x, y)$" is well-defined in **PA**, in any interpretation **M** of **PA**, the well-defined formula "$(Ay)(\sim Q(p, y))$" is a well-defined proposition about an arithmetical relation "$Q(x, y)$" that may, or may not, hold for any, or some, or all values of "$x$" and "$y$" in the domain of **M**. It follows that the arithmetical relation "$Q(x, y)$" is assumed well-defined in every **M**.

However, although "$Q(\underline{k}, \underline{m})$" is equivalent to "$q(k, m)$" in **PA** for any given natural numbers $k, m$ (*where $\underline{k}, \underline{m}$ are the numerals that represent the natural numbers $k, m$ in PA*), the arithmetical relation "$Q(x, y)$" is not equivalent to "$q(x, y)$" in **IA**.

This follows from the fact that whereas "$Q(x, y)$" is a well-formed relation in **IA** that can be expressed entirely in terms of the primitive symbols of **PA**, "$q(x, y)$" is a recursive relation defined in **IA** that cannot be expressed entirely in terms of the primitive symbols of **PA**.

### 6.2. GUS is well-formed and well-defined, but Gödel's reasoning is non-constructive

Hence the proof that "$(E!w)Q(x, w)$" is a well-defined formula of **PA** introduces an element into **PA**, and so into all its interpretations, that is not reflected in the recursive arithmetical relation "$q(x, y)$" of **IA** that "$Q(x, y)$", by the Representation Lemma, is intended to represent formally in **PA**.

This element essentially corresponds to the postulation that the range of values satisfying "$Q(x, y)$", in any interpretation **M** of **PA**, can be used to form a well-defined "set" that completely characterises the defining property "$Q(x, y)$".

I argue elsewhere (*Anand 2001*) that the source of such an obviously non-constructive postulation is the essential use of the *Generalisation Rule of Inference* in the proof that "(E!w)$Q$(x, w)" is a "**PA**-provable" formula.

This rule permits us to infer that the formula "(A$x$)$F$(x)" is "**PA**-provable" if we can establish that the formula "$F$(x)" is "**PA**-provable".

Since the *Generalisation* rule introduces a quantifier (*which ranges over an unspecified domain*) into the formula, I argue that Gödel's proof of the well-definedness of "(E!w)$Q$(x, w)" is non-constructive, even though "(E!w)$Q$(x, w)" is both a well-constructed and a well-defined formula in **PA**.

### 6.3. *Intuitive Arithmetic* **IA** is not a model of **PA**

Since I argue that "$Q$(x, y)" is not a well-defined "relation" of any well-defined set of natural numbers in the above sense, and that the recursive relation "$q$(x, y)" is not "entirely equivalent" to "$Q$(x, y)" in every interpretation of **PA**, I therefore conclude that *Intuitive Arithmetic* **IA** cannot be a model of **PA**.

### 6.4. *Generalisation* as a *Rule of Inference* is intuitionistically questionable

It is in this sense I argue that the Gödel's choice of *Generalisation* as his preferred *Rule of Inference* makes his formal system itself non-constructively Platonistic. The sweeping consequences of *Generalisation* are clearly not obviously rooted in our intuition, and so are intuitionistically questionable.

I argue that *Generalisation* essentially asserts that if "$F$(x)" is "**PA**-provable", then "(A$x$)$F$(x)" is necessarily "formally true" in every interpretation. However, this implicitly assumes that the concept of "intuitive truth" that is native to an interpretation is equivalent to the concept of a "formal truth" in **IA** that is implicit in the definition of "**PA**-provability" based on *Generalisation*.

I argue elsewhere (*Anand 2001, §2.11*) that this implicit assumption is invalid. Thus we can arrive at different assignments of "provability" as a consequence of a definition of "omega-Constructive provability" in a formal system of omega-**PA** where we base the latter on a *Rule of Inference*, such as *Omega-constructivity*, which is inconsistent with *Generalisation*.

### 6.5. A non-constructive interpretation of *Generalisation*

I argue further that if we intend "$(Ax)F(x)$" to interpret as a condensed formula of the assertion "$(Ax)|=_{All\ M} F(x)$" (*in other words of the assertion that "F(x) is 'formally true' for all x in the domain of every interpretation M"*), then Gödel's *Generalisation* is equivalent to the non-formal assertion:

> If "$F(x)$" is "**PA**-provable", then " '$F(x)$ is 'formally true' for all $x$ in the domain of every interpretation **M** of **PA**" is "**PA**-provable".

If we allow that some interpretation **M** of standard **PA** may have a non-countable domain, then *Generalisation* must clearly be viewed as a non-constructive *Rule of Inference*.[11]

Formulas of standard **PA** whose proofs depend on the necessary use of *Generalisation* are thus essentially non-constructive. Their status under interpretation as well-defined "formally true" assertions of any such interpretation **M** must therefore remain in doubt.

---

[11] It essentially postulates, by extrapolation, the "formal truth" of a property, for non-constructive elements in the domain of every interpretation **M** of **PA,** purely on the basis of the logical validity of the property over elements that can be formally represented in **PA**.

## 7. Conclusions

In this paper I argue that, contrary to Gödel's assertion, his argument for the existence of "formally undecidable but intuitively true propositions" in the formal system of standard **PA** (*and other "systems" that are symbolically sufficient to formalise Intuitive Arithmetic recursively*) is not clearly "constructive and intuitionistically unobjectionable".

I argue that the non-constructive element is, in fact, deceptively implicit in the definition of every formal system that admits a *Rule of Inference* such as Gödel's *Generalisation*.

I define a constructive and intuitionistically unobjectionable formal system **PP** that is recursively definable and which represents *Intuitive Arithmetic* (*as symbolised by Peano's Postulates*) more faithfully than Gödel's standard first order **PA**.

I argue that in **PP** (*and similarly definable systems*), Gödel's "undecidable" proposition is a well-formed but ill-defined formal sentence that yields a formal inconsistency similar to that of the "Liar" sentence in ordinary languages.

I also argue that the introduction of "formal truth" and "provability" values to selected propositions of **IA** through **PP** leads to the collapse of the Gödelian argument advanced by J.R.Lucas, Roger Penrose and others.

I argue that the wealth of concepts introduced by Gödel has over-shadowed, and masked, the essentially non-constructive element underlying his reasoning.

I argue that the uncritical acceptance of standard first order **PA** as the most faithful formalisation of *Intuitive Arithmetic*, as expressed by *Peano's Postulates*, must, at some level of consciousness, be subtly impeding the development of more constructive concepts and systems such as **PP** that may be more appropriate for the unfettered development of the emerging area of *Artificial*, and more significantly non-human, *Intelligence*.

Author's e-mail: anandb@vsnl.com